\documentclass[12pt]{article}
\usepackage{mathrsfs}
\usepackage{graphicx,amsmath,amsfonts,amssymb,amscd,amsthm}
\newtheoremstyle{kai}
{3pt} {3pt} {} {} {\bfseries} {.} {.5em} {}
\def\EquationsBySection{\def\theequation
{\thesection.\arabic{equation}}%
\@addtoreset{equation}{section}}

\newcommand\old[1]{}
\newcommand{\pend}{\hfill \thicklines \framebox(6.6,6.6)[l]{}}
\renewenvironment{proof}{\noindent {\it  Proof.} \rm}{\pend}
\newtheorem{theorem}{Theorem}[section]

\newtheorem{remark}{Remark}[section]
\newtheorem{definition}{Definition}[section]
\newtheorem{example}{Example}[section]
\newtheorem{assumption}{Assumption}[section]

\EquationsBySection \makeatother
\usepackage{indentfirst}

\begin{document}
\pagestyle{plain}
\title
{\bf Khasminskii-Type Theorem and LaSalle-Type Theorem for
Stochastic Evolution Delay Equations}

\author{Jianhai Bao$^{a,c}$  Xuerong Mao$^b$, Chenggui Yuan$^{c,}$\thanks{{\it E-mail address:}
C.Yuan@swansea.ac.uk.}\\
\\
${}^a$School of Mathematics, \\Central South
University,Changsha, Hunan 410075, P.R.China \vspace{2mm}\\
 ${}^b$Department of Statistics and Modelling Science,\\
University of Strathclyde, Glasgow  G1 1XH, UK\vspace{2mm}\\
  ${}^c$Department of Mathematics,\\
 University of Wales
Swansea, Swansea SA2 8PP, UK
\\
}

\date{}
\maketitle
\begin{abstract}{\rm In this paper we study the well-known Khasminskii-Type Theorem, i.e. the existence and uniqueness  of solutions of  stochastic
evolution delay equations, under local Lipschitz condition, but without linear
growth condition. We then establish  one
stochastic LaSalle-type theorem for asymptotic stability analysis of
strong solutions.  Moreover, several examples are established to illustrate
the power of our theories.
 }\\

\end{abstract}
\noindent

\section{Introduction}
The study of stochastic evolution delay equations is motivated by
the fact that when one wants to model some evolution phenomena
arising in mechanical, economic, physics, biology, engineering,
etc., some hereditary characteristic such as after-effect, time-lag,
time-delay can appear in the variables ( see, for example, Liu
\cite{liu}, Mohammed \cite{m86} and Wu \cite{w96}). One the other
hand, some of the important and interesting aspects in
existence-and-uniqueness theories and stability analysis for strong
solutions have been greatly developed over the past few years. Here,
we refer to Caraballo et al. \cite{ccm07,clt04,cgr}, Liu
\cite{liu,l98},  Real \cite{r82} and references therein.  For most
of papers mentioned,
 the coefficients of stochastic evolution delay equations
require the global Lipschitz and linear growth conditions
to guarantee  the existence and uniqueness,   and analyze asymptotic
stability for strong solutions. However, there are many stochastic
evolution delay equations which do not  satisfy linear growth
condition, for example:
\begin{equation*}
\begin{cases}
dy(t,x)&=\frac{\partial^2}{\partial
x^2}y(t,x)dt-(y^2(t-\tau,x)-y^3(t,x))dt+y^2(t-\tau,x)dB(t),\ \
t\geq0,\ \
x\in(0,\pi),\\
y(t,x)&=\phi(t,x),\ \ 0\leq x\leq\pi,\ \ t\in[-\tau,0];\ \
y(t,0)=y(t,\pi)=0, \ \ t\geq0,
\end{cases}
\end{equation*}
where $\phi\in C^2([0,\pi]\times[-\tau,0];R)$, $\tau$, positive
constant, and $B(t),t\geq0$, is a real standard Borwnian motion.
Moreover, for such stochastic evolution delay equations, to show
existence-and-uniqueness results and analyze asymptotic stability
for strong solutions, unfortunately, there are not results available
for us to apply. That is, we have no alternative but to put forward
new arguments to overcome the difficulties brought by the nonlinear
growth.

For finite dimensional cases when the drift and diffusion
coefficients of stochastic differential delay equations satisfy
local Lipschitz condition, Mao \cite{m02} established an
existence-and uniqueness theorem of Khasminskii type. Subsequently,
many scholars generalize the classical result to cover
 more general stochastic differential delay equations with Markovian switching and
neutral stochastic differential delay equations, e.g., Mao, Shen and
Yuan \cite{msy08}, Yuan and Glover \cite{yg06} and Yuan and Mao
\cite{ym04}. In particular, it is worth pointing out that
\cite{mr05} by Mao and Rassias gave some Khasminskii-type theorems
for highly nonlinear stochastic differential delay equations and
discussed moment estimations.

On the basis of Khasminskii-type theorems, Mao \cite{m99}
established some stochastic LaSalle-type asymptotic convergence
theorems, and applied to establish sufficient criteria for the
stochastically asymptotic stability of stochastic differential delay
equations. Then, there are extensive literatures which generalize
these stochastic LaSalle-type theorems, see, e.g., Mao, Shen and
Yuan \cite{msy08}, Mao \cite{m01,m02} and Yuan and Mao \cite{ym04}.

However, for stochastic evolution delay equations in infinite
dimensions, as we stated before, in general, the existing
existence-and-uniqueness results and asymptotic stability analysis
for strong solutions are done under the global Lipschitz and linear
growth conditions. Motivated by these papers, we shall intend to
establish one stochastic Khasminskii-type theorem for
existence-and-uniqueness theory and one stochastic LaSalle-type
theorem for asymptotic stability analysis of strong solutions to
stochastic evolution delay equations in infinite dimensions under
local Lipschitz condition, but without linear growth condition. As
we shall see in Section 3 and Section 4, our established theories
have greatly improve some existing results. To the best of our
knowledge to date, there are few literatures concerned with our
problems, therefore, we aim to close a gap.

The contents of this paper will be arranged as follows: In section 2
we collect some preliminaries; In section 3, under local Lipschitz
condition, but without linear growth condition, one Khasminskii-type
theorem is established for stochastic evolution delay equations and
one example is constructed to illustrate the established theory; On
the basis of the established Khasminskii-type theorem, we then
investigate almost surely asymptotic stability for strong solutions,
which is called the LaSalle-type theorem, exponential stability is
also discussed, and
 two examples are provided to explain our theories in the last section.

\section{Preliminaries}
First of all, we introduce the framework in which our analysis is
going to be carried out. Let $V$ be a Banach space and $H, K$ real,
separable Hilbert spaces such that
\begin{equation*}
V\hookrightarrow H\equiv H^{*}\hookrightarrow V^{*},
\end{equation*}
where $V^{*}$ is the dual of $V$ and the injections are continuous ,
dense and compact. We denote by $\|\cdot\|_{*}, \|\cdot\| $ and
$\|\cdot\|_H$ the norms in $V^{*}, V$ and $H$, respectively, by
$\langle\cdot, \cdot\rangle$ the duality product between $V^{*}, V$
and by $\langle\cdot, \cdot\rangle_H$ the scalar product in $H$.
Furthermore, assume that for some $\beta>0$
\begin{equation}\label{eq1}
\beta\| u\|_H\leq\| u\|, \ \ \ \ \ \ \ \ \forall u\in V.
\end{equation}
Assume that $B(t), t\geq0$, is a K-valued Wiener process defined on
a certain probability space $(\Omega, \mathcal {F}, P)$ equipped
with a filtration ${\{\mathcal {F}_t}\}_{t\geq0}$ satisfying the
usual conditions (i.e., it is right continuous and $\mathcal {F}_0$
contains all $\mathcal {P}$-null sets), with covariance operator
$Q\in \mathcal {L}(K)=\mathcal {L}(K,K)$. By \cite[Proposition 4.1,
p.87]{pz92},
\begin{equation*}
E{\langle B(t),x\rangle}_K{\langle B(s),y\rangle}_K=(t\wedge
s){\langle Qx,y\rangle}_K,\ \ \ \ \forall x, y\in K,
\end{equation*}
where $Q$ is a positive, self-adjoint, trace class operator. In
addition, we assume that $e_k, k\in\mathbb{N}$, is an orthonormal
basis of $K$ consisting of eigenvectors of $Q$ with corresponding
eigenvalues $\lambda_k\geq0,k\in\mathbb{N}$, numbered in decreasing
order, and then, according to the representation theorem of
$Q$-Wiener process \cite[Proposition 4.1, p.87]{pz92},
\begin{equation*}
B(t)=\sum\limits_{k=1}^{\infty}\sqrt{\lambda_k}\beta_k(t)e_k,\ \
t\geq0,
\end{equation*}
where $\beta_k(t),k\in\mathbb{N}$ is a sequence of real valued
standard Brownian motions mutually independent on the probability
space $\{\Omega,\mathcal {F},\mathcal {P}\}$. For an operator
$G\in\mathcal {L}(K,H)$, the space of all bounded linear operators
from $K$ into $H$, we denote by $\| G\|_2$ its Hilbert-Schmit norm,
i.e.
\begin{equation*}
\| G\|_2^2=\mbox{trace}(GQG^{*}).
\end{equation*}

In this paper we investigate stochastic evolution delay equation in
the form:
\begin{equation}\label{eq2}
dx(t)=[A(t,x(t))+f(t,x(t),x(t-\tau))]dt+g(t,x(t),x(t-\tau))dB(t),\ \
\ \ \ t\geq0
\end{equation}
with $\tau >0$ and initial datum $x(\theta)=\psi(\theta)\in
C^b_{\mathcal {F}_0}([-\tau,0];V)\cap C^b_{\mathcal
{F}_0}([-\tau,0];H)$, the family of all $\mathcal {F}_0$-measurable
bounded $C([-\tau,0];V)\cap C([-\tau,0];H)$-valued random variables.
\begin{assumption}\label{as1}
For any $T>0$ let $A(t,\cdot):V\rightarrow V^*$ be a family of
(nonlinear) operators on $t\in[0,T]$ satisfying $A(t,0)=0$ and
$p\geq2$:
\item[\textmd{(A.1)}] (Monotonicity and Coercivity) $\exists
\alpha>0, \lambda \in R$ such that
\begin{equation*}
2\langle A(t,x)-A(t,y),x-y\rangle\leq-\alpha \| x-y\|^p +\lambda\|
x-y\|_H^2,\ \ \ \ \ \ \forall x,y\in V;
\end{equation*}
\item[\textmd{(A.2)}] (Measurability) $\forall x\in V$, the map $t\in[0,T]\rightarrow A(t,x)\in V^{*}$
 is Lebesgue measurable;
\item[\textmd{(A.3)}](Hemicontinuity) The map
\begin{equation*}
\theta\in R\rightarrow\langle A(t,x+\theta y), z\rangle\in R
\end{equation*}
is continuous for arbitrary $x,y,z\in V$ and $t\in[0,T]$;
\item[\textmd{(A.4)}] (Boundedness) $\exists$ $\gamma>0$ such that
for $t\in[0,T]$
\begin{equation*}
\|A(t,x)\|_{*}\leq \gamma\| x\|^{p-1}, \ \ \ \ \ \ \forall x\in V.
\end{equation*}
\end{assumption}

\begin{assumption}\label{as2}
 Let $f(t,\cdot,\cdot):H\times H\rightarrow H$ and $g(t,\cdot,\cdot):H\times H\rightarrow \mathcal {L}(K,H)$  be the
 families of nonlinear operators defined for $t\in[0,T]$ and satisfy:
\begin{enumerate}
\item[\textmd{(B.1)}] (Measurability) For any $x,y\in H$, the maps
\begin{equation*}
\forall t\in[0,T]\rightarrow f(t,x,y)\in H\ \ \ \ \mbox{and}\ \ \
g(t,x,y)\in \mathcal {L}(K,H)
\end{equation*}
are Lebesgue-measurable, respectively.

\item[\textmd{(B.2)}] (Boundedness)
\begin{equation}\label{0}
M:=\sup\limits_{t\geq0}\{\|f(t,0,0)\|_H\vee\|g(t,0,0)\|_2\}<\infty.
\end{equation}
\end{enumerate}

\end{assumption}

\begin{assumption}\label{as3}
(Local Lipschitz Condition) For each $h>0$, $\exists L_h>0$ such
that
\begin{equation}\label{eq0}
\|f(t,x_1,y_1)-f(t,x_2,y_2)\|_H\vee\|g(t,x_1,y_1)-g(t,x_2,y_2)\|_2\leq
L_h(\|x_1-y_1\|_H+\|x_2-y_2\|_H)
\end{equation}
for $t\in[0,T]$ and $x_1,x_2,y_1,y_2\in H$ with
$\|x_1\|_H\vee\|x_2\|_H\vee\|y_1\|_H\vee\|y_2\|_H\leq h$;
\end{assumption}

Let $I^p([-\tau, T];V)$ denotes the space of all $V$-valued
processes $x(t)$, which are $\mathcal {F}_t$-measurable from
$[-\tau,T]\times\Omega$ to $V$ and satisfy $
E\int_{-\tau}^T\|x(t)\|^pdt<\infty, $ and $L^2(\Omega;
C([-\tau,T];H))$ denote $H$-valued processes from
$[-\tau,T]\times\Omega$ to $H$ such that $E\sup_{-\tau\le t\le
T}\|x(t)\|_H^2<\infty.$ Let's recall the definition of strong
solutions.
\begin{definition}
For any initial datum $x(\theta)=\psi(\theta)\in C^b_{\mathcal
{F}_0}([-\tau,0];V)\cap C^b_{\mathcal {F}_0}([-\tau,0];H)$, a
stochastic process $x(t), t\in[0,T]$, is said to be a strong
solution of \eqref{eq2} if the following conditions are satisfied:
\begin{enumerate}
\item[\textmd{(a)}]  $x(t)\in I^p([-\tau ,T];V)\bigcap
L^2(\Omega; C([-\tau,T];H))$;

\item[\textmd{(b)}] The following equality holds in $V^*$ almost
surely for $t\in[0,T]$
\begin{equation*}
x(t)=x(0)+\int_0^t[A(s,x(
s))+f(s,x(s),x(s-\tau))]ds+\int_0^tg(s,x(s),x(s-\tau))dB(s)
\end{equation*}
with initial condition $x(\theta)=\psi(\theta)\in C^b_{\mathcal
{F}_0}([-\tau,0];V)\cap C^b_{\mathcal {F}_0}([-\tau,0];H)$.
\end{enumerate}
If $T$ is replaced by $\infty$, $x(t),t\geq0$, is called a global
strong solution of \eqref{eq2}.

\end{definition}

In what follows, we shall also need the following global Lipschitz
condition.
\begin{assumption}\label{as4}
There exists a constant $L>0$ such that, for $t\in[0,T]$ and
$x_1,x_2,y_1,y_2\in H$,
\begin{equation}\label{eq4}
\|f(t,x_1,y_1)-f(t,x_2,y_2)\|_H+\|g(t,x_1,y_1)-g(t,x_2,y_2)\|_2\leq
L(\|x_1-y_1\|_H+\|x_2-y_2\|_H).
\end{equation}
\end{assumption}

Under the global Lipschitz condition \eqref{eq4}, the following
existence-and-uniqueness result can be found in \cite[Theorem
3.1]{cgr}.

\begin{theorem}\label{global}
Assume that Assumption \ref{as1}, Assumption \ref{as2} and
Assumption \ref{as4} hold. Then, for each $\psi\in C^b_{\mathcal
{F}_0}([-\tau,0];V)\cap C^b_{\mathcal {F}_0}([-\tau,0];H)$ there
exists a unique strong solution of \eqref{eq2} in
$I^p([-\tau,T];V)\bigcap L^2(\Omega;C([-\tau,T];H))$.
\end{theorem}

\begin{remark}
In general, Assumption \ref{as1}-\ref{as3} will only guarantee a
unique maximal local strong solution to \eqref{eq2} for any given
initial data $\psi$. However, the additional conditions imposed in
one of our main results, Theorem \ref{existence}, will guarantee
that this maximal local strong solution is in fact a unique global
strong one.
\end{remark}

Now we recall the It\^o formula which will play key role in what
follows. Let $R_+$ be non-negative real number, and
$C^{1,2}(R_+\times H;R_+)$ denote the space of all real valued
non-negative functions $U$ on $R_+\times H$ with properties:
\begin{enumerate}
\item[\textmd{(i)}] $U(t,x)$ is once differential in $t$ and twice
(Fr\'echet) differentiable in $x$;
\item[\textmd{(i)}] $U_x(t,x)$ and $U_{xx}(t,x)$ are both continuous
in $H$ and $L(H)$, respectively.
\end{enumerate}
\begin{theorem}\label{yito}
Suppose $U\in C^{1,2}(R_+\times H;R_+)$ and $x(t),t\geq0$, is a
strong solution to \eqref{eq2}, then
\begin{equation*}
\begin{split} U(t,x(t))&=U(0,\psi(0))+\int_0^t\mathcal
{L}U(s,x(s),x(s-\tau))ds+\int_0^t\langle
U_x(s,x(s)),g(s,x(s),x(s-\tau))dB(s)\rangle_H,\\
\end{split}
\end{equation*}
where $\mathcal {L}$ is the associated diffusion operator defined
by, for any $t\geq0$ and $x,y\in V$ ,
\begin{equation*}
\begin{split}
\mathcal {L}U(t,x,y)&=\frac{\partial U(t,x)}{\partial t}+\langle
A(t,x)+f(t,x,y),U_x(t,x)
\rangle+\frac{1}{2}\mbox{trace}(U_{xx}(t,x)g(t,x,y)Qg^*(t,x,y)).
\end{split}
\end{equation*}
\end{theorem}

We will also need the following useful semimartingale convergence
theorem established by Lipster and Shiryayev \cite[Theorem 7,
p.139]{ls89}.

\begin{theorem}\label{Lipster}
Let $A_1(t)$ and $A_2(t)$ be two continuous adapted increasing
processes on $t\geq0$ with $A_1(0)=A_2(0)=0$ a.s. Let $M(t)$ be a
real-valued continuous local martingale with $M(0)=0$ a.s. Let
$\zeta$ be a nonnegative $\mathcal {F}_0$-measurable random variable
such that $E\zeta<\infty$. Define for $t\geq0$
\begin{equation*}
X(t):=\zeta + A_1(t)-A_2(t)+M(t).
\end{equation*}
If $X(t)$ is nonnegative, then a.s.
\begin{equation*}
\left\{\lim\limits_{t\rightarrow\infty}A_1(t)<\infty\right\}\subset\left\{\lim\limits_{t\rightarrow\infty}X(t)<\infty\right\}\bigcap
\left\{\lim\limits_{t\rightarrow\infty}A_2(t)<\infty\right\},
\end{equation*}
where $C\subset D$ a.s. means $P(C\bigcap D^c)=0$. In particular, if
$\lim\limits_{t\rightarrow\infty} A_1(t)<\infty$ a.s., then, with
probability one,
\begin{equation*}
\lim\limits_{t\rightarrow\infty} X(t)<\infty,\ \ \ \ \ \ \ \ \
\lim\limits_{t\rightarrow\infty} A_2(t)<\infty
\end{equation*}
and
\begin{equation*}
-\infty<\lim\limits_{t\rightarrow\infty} M(t)<\infty.
\end{equation*}
That is, all of the three processes $X(t), A_2(t)$ and $M(t)$
converge to finite random variables.
\end{theorem}

\section{Khasminskii-Type Theorem}

In this section, {\it under local Lipschitz condition, but without
linear growth condition}, we shall establish one Khasminskii-type
theorem for existence-and uniqueness theory for stochastic evolution
delay equations in infinite dimensions.

\begin{theorem}\label{existence}
Let Assumption \ref{as1}-\ref{as3} hold. Assume further that there
are functions $U\in C^{1, 2}(R_+\times H;R_+)$, $W\in C(R_+\times
H;R_+)$, and positive constants $\lambda_1$ and $\lambda_2$ such
that
\begin{equation}\label{eq7}
\mathcal {L}U(t,x,y)\leq
\lambda_1[1+U(t,x)+U(t-\tau,y)+W(t-\tau,y)]-\lambda_2W(t,x), \ \ \
(t,x,y)\in R_+\times V\times V,
\end{equation}
and
\begin{equation}\label{eq8}
\lim\limits_{\|x\|_H\rightarrow\infty}\inf\limits_{0\leq
t<\infty}U(t,x)=\infty,\ \ \ \ \ x\in V.
\end{equation}
Then, for any initial data $x(\theta)=\psi(\theta)\in C^b_{\mathcal
{F}_0}([-\tau,0];V)\cap C^b_{\mathcal {F}_0}([-\tau,0];H)$,
\eqref{eq2} admits a unique global solution.
\end{theorem}

\begin{proof}
For any integer $k\geq b$, bound of $\psi$, and $x,y\in H$, define
\begin{equation*}
f_k(t,x,y)=f\left(t,\frac{\|x\|_H\wedge
k}{\|x\|_H}x,\frac{\|y\|_H\wedge k}{\|y\|_H}y\right),\ \ \ \ \
g_k(t,x,y)=g\left(t,\frac{\|x\|_H\wedge
k}{\|x\|_H}x,\frac{\|y\|_H\wedge k}{\|y\|_H}y\right),
\end{equation*}
where we set $(\|x\|_H\wedge k/\|x\|_H)x=0$ when $x=0$. Then, by
\eqref{eq0} and \eqref{0}, for any $x,y\in H$ and $t\geq0$ we
observe that $f_k(t,x,y)$ and $g_k(t,x,y)$ satisfy the global
Lipschitz condition and linear growth condition. Hence, there exists
by Theorem \ref{global} a unique global solution $x_k(t)$ on
$[-\tau,\infty)$ to the following stochastic evolution delay
equation
\begin{equation*}
dx_k(t)=[A(t,x_k(t))+f_k(t,x_k(t),x_k(t-\tau))]dt+g_k(t,x_k(t),x_k(t-\tau))dB(t)
\end{equation*}
with initial data $x(\theta)=\psi(\theta)\in C^b_{\mathcal
{F}_0}([-\tau,0];V)\cap C^b_{\mathcal {F}_0}([-\tau,0];H)$. Define
the stopping time
\begin{equation*}
\sigma_k=\inf\{t\geq0:\|x_k(t)\|_H\geq k\},
\end{equation*}
where we set $\inf\emptyset=\infty$ as usual. Clearly, for any
$s\leq\sigma_k$, $\|x_k(s)\|_H\wedge\|x_k(s-\tau)\|_H\leq k$. Then,
recalling the definition of $f_k$ and $g_k$, it is easy to see that,
for any $0\leq s\leq \sigma_k$,
\begin{equation*}
f_{k+1}(s,x_k(s),x_k(s-\tau))=f_k(s,x_k(s),x_k(s-\tau))=f(s,x_{k+1}(s),x_{k+1}(s-\tau))
\end{equation*}
and
\begin{equation*}
g_{k+1}(s,x_k(s),x_k(s-\tau))=g_k(s,x_k(s),x_k(s-\tau))=g(s,x_{k+1}(s),x_{k+1}(s-\tau)).
\end{equation*}
Consequently,
\begin{equation*}
\begin{split}
x_k&(t\wedge\sigma_k)\\
&=\psi(0)+\int_0^{t\wedge\sigma_k}[A(s,x_k(
s))+f_k(s,x_k(s),x_k(s-\tau))]ds+\int_0^{t\wedge\sigma_k}g_k(s,x_k(s),x_k(s-\tau))dB(s)\\
&= \psi(0)+\int_0^{t\wedge\sigma_k}[A(s,x_k(
s))+f_{k+1}(s,x_k(s),x_k(s-\tau))]ds+\int_0^{t\wedge\sigma_k}g_{k+1}(s,x_k(s),x_k(s-\tau))dB(s),
\end{split}
\end{equation*}
which immediately gives
\begin{equation*}
x_{k+1}(t)=x_k(t),\ \ \ \ \ \ \ \ 0\leq t\leq \sigma_k.
\end{equation*}
This further implies that $\sigma_k$ is increasing in k. So we can
define $\sigma= \lim\limits_{k\rightarrow\infty}\sigma_k$ . The
property above also enables us to define $x(t)$ for $t\in[-\tau,
\sigma)$ as follows
\begin{equation*}
x(t)=x_k(t),\ \ \ \ \ \ \ \ -\tau\leq t\leq \sigma_k.
\end{equation*}
It is clear that $x(t)$ is a unique solution to \eqref{eq2} for
$t\in [-\tau, \sigma_k)$. To complete the proof, we only need to
show that $P(\sigma=\infty)=1$. Indeed, to show the desired
assertion we compute by the It\^o formula and \eqref{eq7} that for
any $t\in[0,\tau]$
\begin{equation}\label{eq11}
\begin{split}
EU(t\wedge\sigma_k,x(t\wedge\sigma_k))&=EU(0,\psi(0))+E\int_0^{t\wedge\sigma_k}\mathcal
{L}U(s,x(s),x(s-\tau))ds\\
&\leq
EU(0,\psi(0))+E\int_0^{\tau}\lambda_1[1+U(s-\tau,x(s-\tau))+W(s-\tau,x(s-\tau))]ds\\
&+\lambda_1E\int_0^{t\wedge\sigma_k}U(s,x(s))ds-\lambda_2E\int_0^{t\wedge\sigma_k}W(s,x(s))ds\\
&\leq
C_1+\lambda_1E\int_0^{t}U(s\wedge\sigma_k,x(s\wedge\sigma_k))ds-\lambda_2E\int_0^{t\wedge\sigma_k}W(s,x(s))ds,
\end{split}
\end{equation}
where
\begin{equation*}
C_1=EU(0,\psi(0))+E\int_{-\tau}^0\lambda_1[1+U(s,\psi(s))+W(s,\psi(s))]ds.
\end{equation*}
Therefore, for $t\in[0,\tau]$ the Gronwall inequality yields
\begin{equation}\label{eq9}
EU(t\wedge\sigma_k,x(t\wedge\sigma_k))\leq C_1e^{\lambda_1\tau},
\end{equation}
and, in addition to the definition of $\sigma_k$,
\begin{equation*}
P(\sigma_k\leq \tau)\leq
\frac{C_1e^{\lambda_1\tau}}{\inf_{t\geq0,\|x\|_H\geq k}U(t,x)}.
\end{equation*}
Letting $k\rightarrow\infty$ and then observing \eqref{eq8}, we
obtain
\begin{equation*}
P(\sigma\leq \tau)=0.
\end{equation*}
Namely,
\begin{equation}\label{eq12}
P(\sigma>\tau)=1.
\end{equation}
Hence, in \eqref{eq9}, letting $k\rightarrow\infty$ gives, for any
$t\in[0,\tau]$ ,
\begin{equation}\label{eq10}
EU(t,x(t))\leq C_1e^{\lambda_1\tau}.
\end{equation}
However, by \eqref{eq11}, for any $t\in[0,\tau]$,
\begin{equation*}
\lambda_2E\int_0^{t\wedge\sigma_k}W(s,x(s))ds\leq
C_1+\lambda_1E\int_0^tU(s\wedge\sigma_k,x(s\wedge\sigma_k))ds.
\end{equation*}
So, letting $k\rightarrow\infty$, together with \eqref{eq12} and
\eqref{eq10}, it could be deduced that
\begin{equation}\label{eq13}
E\int_0^{\tau}W(x(s))ds\leq \frac{C_1+\lambda_1\tau
C_1e^{\lambda_1\tau}}{\lambda_2}<\infty.
\end{equation}
In the same manner as \eqref{eq11} was done, for any $t\in[0,2\tau]$
\begin{equation}\label{eq14}
EU(t\wedge\sigma_k,x(t\wedge\sigma_k))\leq C_2
+\lambda_1E\int_0^{t}U(s\wedge\sigma_k,x(s\wedge\sigma_k))ds-\lambda_2E\int_0^{t\wedge\sigma_k}W(s,x(s))ds,
\end{equation}
where
\begin{equation}
\begin{split}
C_2&=EU(0,\psi(0))+E\int_0^{2\tau}\lambda_1[1+U(s-\tau,x(s-\tau))+W(s-\tau,x(s-\tau))]ds\\
&=EU(0,\psi(0))+E\int_{-\tau}^0\lambda_1[1+U(s,x(s))+W(s,x(s))]ds\\
&+E\int_0^{\tau}\lambda_1[1+U(s,x(s))+W(s,x(s))]ds\\
&<\infty
\end{split}
\end{equation}
since \eqref{eq10} and \eqref{eq13} hold for any $t\in[0,\tau]$. We
then have by Gronwall's inequality that, for any $t\in[0,2\tau]$,
\begin{equation}\label{eq15}
EU(t\wedge\sigma_k,x(t\wedge\sigma_k))\leq C_2e^{2\lambda_1\tau}.
\end{equation}
Therefore, taking into account \eqref{eq8},
\begin{equation*}
P(\sigma\leq 2\tau)=0,
\end{equation*}
that is
\begin{equation*}
P(\sigma> 2\tau)=1.
\end{equation*}
Next, by letting $k\rightarrow\infty$ in \eqref{eq14} and
\eqref{eq15}, respectively, we derive that, for any $t\in[0,2\tau]$,
\begin{equation*}
EU(t,x(t))\leq C_2e^{2\lambda_1\tau}\ \ \ \mbox{ and }  \ \ \
E\int_0^{2\tau}W(x(s))ds<\frac{C_2+2\lambda_1\tau
C_2e^{2\lambda_1\tau}}{\lambda_2}<\infty.
\end{equation*}
By induction, for any integer $k\geq1$ it follows easily that
\begin{equation*}
EU(t,x(t))\leq C_ke^{k\lambda_1\tau}
\end{equation*}
whenever $t\in[0,k\tau]$ and
\begin{equation*}
E\int_0^{k\tau}W(s,x(s))ds<\frac{C_k+k\lambda_1\tau
C_ke^{k\lambda_1\tau}}{\lambda_2}.
\end{equation*}
So, we can conclude by \eqref{eq8} that
\begin{equation*}
P(\sigma<\infty)=0,
\end{equation*}
and then \eqref{eq2} admits a unique global strong solution on
$t\geq0$.
\end{proof}

 Now we shall use   Theorem
\ref{existence} to analyze the example which appeared  the introduction section.

\begin{example}
Consider the following semilinear stochastic partial differential
equation:
\begin{equation}\label{eq16}
\begin{cases}
dy(t,x)&=\frac{\partial^2}{\partial
x^2}y(t,x)dt+(y^2(t-\tau,x)-y^3(t,x))dt+y^2(t-\tau,x)dB(t),\ \
t\geq0,\ \
x\in(0,\pi),\\
y(t,x)&=\phi(t,x),\ \ 0\leq x\leq\pi,\ \ t\in[-\tau,0];\ \
y(t,0)=y(t,\pi)=0, \ \ t\geq0,
\end{cases}
\end{equation}
where $\phi\in C^2([0,\pi]\times[-\tau,0];R)$, $\tau$, positive
constant, and $B(t),t\geq0$, is a real standard Borwnian motion.
Take $H=L^2([0,\pi]), V=H_0^1([0,\pi]), V^*=H^{-1}([0,\pi]), K=R,
A(t,u)=\frac{\partial^2}{\partial x^2}u(x),f(t,u,v)=v^2(x)-u(x)$ and
$g(t,u,v)=v^2(x), u,v\in V$. Furthermore, the norms in $H$ and $V$
are defined as
$\|\xi\|_H=\left(\int_0^{\pi}\xi^2(s)ds\right)^{\frac{1}{2}}$ for
$\xi\in H$ and $\|\xi\|=\left(\int_0^{\pi}\left(\frac{\partial
\xi}{\partial s}\right)^2 ds\right)^{\frac{1}{2}}$ for $\xi\in V$.
Then, clearly, in \eqref{eq1}, we can take $\beta=1$. Setting
$U(t,x)=\|x\|_H^2$ and recalling the definition of diffusion
operator $\mathcal {L}U$, it follows easily that
\begin{equation*}
\begin{split}
\mathcal {L}U(t,x,y)&=2\langle A(t,x),x\rangle+2\langle
f(t,x,y),x\rangle_H+\|g(t,x,y)\|_H^2\\
&=2\langle A(t,x),x\rangle+2\langle y^2-x^3,x\rangle_H+\|
y^2\|_H^2\\
& \leq-2\|x\|^2+2\|x\|_H\|y\|^2_H-2\|x\|^4_H+\|y\|^4_H\\
&\leq\|x\|^2_H-2\|x\|^4_H+\frac{4}{3}\|y\|^4_H.
\end{split}
\end{equation*}
Hence, by Theorem \ref{existence}, setting $\lambda_1=1$ and
$\lambda_2=\frac{4}{3}$ and $W(t,x)=\|x\|_H^4$, we immediately
conclude that \eqref{eq16} admits a global solution on $t\geq0$.
\end{example}

\begin{remark}
Since $f$ and $g$ do not satisfy the linear growth condition, then
\cite[Theorem 3.1]{clt04} certainly cannot apply to the above
example. However, by our established theory we can deduce that
\eqref{eq16} has a unique global solution on $t\geq0$. Therefore,
Theorem \ref{existence} covers many highly nonlinear stochastic
evolution delay equations.
\end{remark}

\section{LaSalle-Type Theorem}
On the basis of the established Khasminskii-type theorem, in what
follows we shall analyze the asymptotic stability properties of
strong solutions under some special conditions of Khasminskii type by using Lyapunov method.
As we know, the Lyapunov method has been developed and applied by many authors during the past century.
In 1968, (see \cite{l68}, Hale and Lunel \cite{hl93} and the reference therein),  LaSalle used the
the Lyapunov method to locate limit sets for ordinary nonautonomous systems, which is one of the
 important developments in this direction, and the theorem is called
 LaSalle theorem. After thirty 30 years, Mao \cite{m99} established a stochastic version of
 the LaSalle theorem for stochastic differential equations in finite dimensional space.
 In this section, we shall extend the LaSalle theorem to the strong solutions of stochastic evolution delay equations.
  We shall see there are a lot of difficulties to overcome from
 finite dimensional cases to infinite dimensional cases.

Let $L^1(R_+;R_+)$ denote the family of all functions $\xi: R_+
\rightarrow R_+$ such that $\int_0^{\infty}\xi(s)ds<\infty$.
\begin{theorem}\label{Lasalle}
Let Assumption \ref{as1}-\ref{as3} hold. Assume that there are
functions $U\in C^{1, 2}(R_+\times H;R_+)$, $\gamma\in L^1(R_+;R_+)$
and $w_1,w_2\in C(H;R_+)$ such that
\begin{equation}\label{4}
\mathcal {L}U(t,x,y)\leq \gamma(t)-w_1(x)+w_2(y), \ \ \ \ \ \ \
\forall (t,x,y)\in R_+\times V\times V,
\end{equation}
\begin{equation}\label{5}
w_1(0)=w_2(0)=0,\ \ \ \ \ \ \ \ w_1(x)>w_2(x),
\end{equation}
moreover
\begin{equation}\label{6}
\lim\limits_{\|x\|_H\rightarrow\infty}\inf\limits_{0\leq
t<\infty}U(t,x)=\infty \ \mbox{ and }\ \
\lim\limits_{\|x\|\rightarrow\infty}\inf\limits_{0\leq
t<\infty}U(t,x)=\infty,\ \ \ \ \  x\in V.
\end{equation}
Then the solution $x(t)$ of \eqref{eq2} satisfies
\begin{equation}\label{36}
\lim\limits_{t\rightarrow\infty}\sup U(t,x(t))<\infty\ \ \ \ \ \
a.s.,
\end{equation}
and
\begin{equation}\label{37}
\lim\limits_{t\rightarrow\infty}
w(x(t))=0 \ \ \ \ \ \ a.s.,
\end{equation}
where $w=w_1-w_2$, and moreover
\begin{equation}\label{38}
P\left(\lim\limits_{t\rightarrow\infty}\|x(t)\|_H=0\right)=1,
\end{equation}
that is, the solution of \eqref{eq2} is almost surely asymptotically
stable.
\end{theorem}

\begin{proof}
First of all,  by Theorem \ref{existence}  it is easy to see
\eqref{eq2} has a unique global solution for $t\geq0$  under the
conditions of Theorem \ref{Lasalle}. Applying the It\^o formula to
$V(t,x)$ and solution $x(t), t\geq0,$ of \eqref{eq2}, we derive that
\begin{equation*}
U(t,x(t))=U(0,\psi(0))+\int_0^t\mathcal
{L}U(s,x(s),x(s-\tau))ds+\int_0^t\langle
U_x(s,x(s)),g(s,x(s),x(s-\tau))dB(s)\rangle_H.
\end{equation*}
This, together with \eqref{4},  implies that
\begin{equation}\label{7}
\begin{split}
U(t,x(t))&\leq
U(0,\psi(0))+\int_0^t[\gamma(s)-w_1(x(s))+w_2(x(s-\tau))]ds\\
&+\int_0^t\langle U_x(s,x(s)),g(s,x(s),x(s-\tau))dB(s)\rangle_H.\\
&=U(0,\psi(0))+\int_0^t\gamma(s)ds+\int_{-\tau}^0w_2(\psi(s))ds\\
&-\int_0^t[w_1(x(s))-w_2(x(s))]ds +\int_0^t\langle
U_x(s,x(s)),g(s,x(s),x(s-\tau))dB(s)\rangle_H.
\end{split}
\end{equation}
Then, by \eqref{5} and Theorem \ref{Lipster}  we obtain
\begin{equation}\label{8}
\lim\limits_{t\rightarrow\infty}\sup U(t,x(t))<\infty\ \ \ \ \ \
\mbox{ a.s}.
\end{equation}
Taking expectations on both sides of \eqref{7} and then letting
$t\rightarrow\infty,$ one derives that
\begin{equation}\label{9}
E\int_0^{\infty}(w_1(x(s))-w_2(x(s)))ds<\infty,
\end{equation}
which certainly implies
\begin{equation}\label{10}
\int_0^{\infty}(w_1(x(s))-w_2(x(s)))ds<\infty,\ \ \ \ \ \
\mbox{a.s}.
\end{equation}
Clearly, $w\in C(H;R_+)$. It is straightforward to see from
\eqref{10} that
\begin{equation}\label{11}
\lim\limits_{t\rightarrow\infty}\inf w(x(t))=0 \ \ \ \ \ \ \mbox{
a.s}.
\end{equation}
In what follows, we intend to claim that
\begin{equation}\label{12}
\lim\limits_{t\rightarrow\infty} w(x(t))=0 \ \ \ \ \ \ \mbox{a.s}.
\end{equation}
By contradiction, if \eqref{12} is false, then
\begin{equation}\label{13}
P\left\{\lim\limits_{t\rightarrow\infty} \sup w(x(t))>0\right\} >0.
\end{equation}
Hence, there is a number $\epsilon>0$ such that
\begin{equation}\label{14}
P(\Omega_1)\geq3\epsilon,
\end{equation}
where
\begin{equation*}
\Omega_1=\left\{\lim\limits_{t\rightarrow\infty} \sup
w(x(t))>2\epsilon\right\}.
\end{equation*}
It is easy to observe from \eqref{8} and the continuity of both the
solution $x(t)$ and the function $U(t,x)$ that
\begin{equation*}
\sup U(t,x(t))<\infty\ \ \ \ \ \ \mbox{ a.s}.
\end{equation*}
Define $\rho:R_+\rightarrow R_+$ by
\begin{equation*}
\rho(r)=\inf\limits_{\|x\|\geq r,\  0\leq t<\infty}U(t,x).
\end{equation*}
Clearly, $\rho(\|x(t)\|)\leq U(t,x(t))$ so
\begin{equation*}
\sup\limits_{0\leq t<\infty}\rho(\|x(t)\|)\leq\sup\limits_{0\leq
t<\infty}U(t,x(t))<\infty\ \ \ \ \ \ \mbox{ a.s}.
\end{equation*}
While by \eqref{6}
\begin{equation*}
\rho(r)=\infty.
\end{equation*}
We therefore must have
\begin{equation}\label{15}
\sup\limits_{0\leq t<\infty}\|x(t)\|<\infty\ \ \ \ \ \ a.s.
\end{equation}
Recalling the boundedness of the initial data we can then find a
positive number $h$, which depends on $\epsilon$, sufficiently large
for $\|\psi(\theta)\|<h$ for all $-\tau\leq\theta \leq 0$ almost
surely while
\begin{equation}\label{16}
P\left(\Omega_2\right)\geq1-\epsilon,
\end{equation}
where
\begin{equation*}
\Omega_2=\left\{\sup\limits_{-\tau\leq t<\infty}\|x(t)\|<h\right\}.
\end{equation*}
It is easy to see from \eqref{14} and \eqref{16} that
\begin{equation}\label{17}
P\left(\Omega_1\cap\Omega_2\right)\geq2\epsilon.
\end{equation}
Let us now define a sequence of stopping times,
\begin{equation*}
\begin{split}
\tau_h&=\inf\{t\geq0:\|x(t)\|\geq h\},\\
\sigma_1&=\inf\{t\geq0:w(x(t))\geq2\epsilon\},\\
\sigma_{2k}&=\inf\{t\geq\sigma_{2k-1}:w(x(t))<\epsilon\},\ \ \ \ k=1,2,\cdots\\
\sigma_{2k+1}&=\inf\{t\geq\sigma_{2k}:w(x(t))\geq2\epsilon\},\ \ \ \ k=1,2,\cdots,\\
\end{split}
\end{equation*}
where throughout this paper we set $\inf\emptyset=\infty$. Note from
\eqref{11} and the definitions of $\Omega_1$ and $\Omega_2$ that
\begin{equation}\label{27}
\tau_h=\infty, \ \ \ \sigma_k<\infty, \ \ \forall k\geq1
\end{equation}
whenever $\omega\in\Omega_1\bigcap\Omega_2$. By \eqref{9}, we
compute
\begin{equation}\label{31}
\begin{split}
\infty&>E\int_0^{\infty}w(t,x(t))dt\\
&\geq\sum\limits_{k=1}^{\infty}E\left[I_{\{\sigma_{2k-1}<\infty,\sigma_{2k}<\infty,\tau_h=\infty\}}\int_{\sigma_{2k-1}}^{\sigma_{2k}}w(t,x(t))dt\right]\\
&\geq\epsilon\sum\limits_{k=1}^{\infty}E[I_{\{\sigma_{2k-1}<\infty,\tau_h=\infty\}}(\sigma_{2k}-\sigma_{2k-1})],
\end{split}
\end{equation}
where $I_A$ is the indicator function of set $A$ and we have noted
from \eqref{11} that ¦Ò$\sigma_{2k}<\infty$ whenever
$\sigma_{2k-1}<\infty$. On the other hand, by It\^o's formula and
$(A.1)$
\begin{equation}\label{18}
\begin{split}
E&\left[I_{\{\tau_h\wedge\sigma_{2k-1}<\infty\}}\sup\limits_{0\leq
t\leq T}
\|x(\tau_h\wedge(\sigma_{2k-1}+t))-x(\tau_h\wedge\sigma_{2k-1})\|_H^2\right]\\
&+\alpha
E\left[I_{\{\tau_h\wedge\sigma_{2k-1}<\infty\}}\sup\limits_{0\leq
t\leq
T}\int_{\tau_h\wedge\sigma_{2k-1}}^{\tau_h\wedge(\sigma_{2k-1}+t)}\|x(s)\|^2ds\right]\\
&\leq|\lambda|E\left[I_{\{\tau_h\wedge\sigma_{2k-1}<\infty\}}\int_{\tau_h\wedge\sigma_{2k-1}}^{\tau_h\wedge(\sigma_{2k-1}+T)}\|x(s)\|_H^2ds\right]\\
&+
2E\left[I_{\{\tau_h\wedge\sigma_{2k-1}<\infty\}}\int_{\tau_h\wedge\sigma_{2k-1}}^{\tau_h\wedge(\sigma_{2k-1}+T)}
|\langle A(s,x(s)),x(\tau_h\wedge\sigma_{2k-1})\rangle |ds\right]\\
&+2E\left[I_{\{\tau_h\wedge\sigma_{2k-1}<\infty\}}\int_{\tau_h\wedge\sigma_{2k-1}}^{\tau_h\wedge(\sigma_{2k-1}+T)}
|\langle f(s,x(s),x(s-\tau)),x(s)-x(\tau_h\wedge\sigma_{2k-1})\rangle_H| ds\right]\\
&+E\left[I_{\{\tau_h\wedge\sigma_{2k-1}<\infty\}}\int_{\tau_h\wedge\sigma_{2k-1}}^{\tau_h\wedge(\sigma_{2k-1}+T)}\|g(s,x(s),x(s-\tau))\|_2^2ds\right]\\
&+2E\left[I_{\{\tau_h\wedge\sigma_{2k-1}<\infty\}}\sup\limits_{0\leq
t\leq
T}\left|\int_{\tau_h\wedge\sigma_{2k-1}}^{\tau_h\wedge(\sigma_{2k-1}+t)}\langle
x(s)-x(\tau_h\wedge\sigma_{2k-1}),
g(s,x(s),x(s-\tau))dB(s)\rangle_H\right|\right].\\
\end{split}
\end{equation}
Obviously, it follows easily from \eqref{eq1} that
\begin{equation}\label{19}
\begin{split}
|\lambda|E\left[I_{\{\tau_h\wedge\sigma_{2k-1}<\infty\}}
\int_{\tau_h\wedge\sigma_{2k-1}}^{\tau_h\wedge(\sigma_{2k-1}+T)}\|x(s)\|_H^2ds\right]
\leq\frac{|\lambda|Th^2}{\beta^2}:=C_1T^{\frac{1}{2}}.
\end{split}
\end{equation}
Now compute by $(A.4)$ that
\begin{equation}\label{20}
\begin{split}
2&E\left[I_{\{\tau_h\wedge\sigma_{2k-1}<\infty\}}\int_{\tau_h\wedge\sigma_{2k-1}}^{\tau_h\wedge(\sigma_{2k-1}+T)}
|\langle A(s,x(s)),x(\tau_h\wedge\sigma_{2k-1})\rangle |ds\right]\\
&\leq2E\left[I_{\{\tau_h\wedge\sigma_{2k-1}<\infty\}}\int_{\tau_h\wedge\sigma_{2k-1}}^{\tau_h\wedge(\sigma_{2k-1}+T)}
\|A(s,x(s))\|_*\|x(\tau_h\wedge\sigma_{2k-1})\|ds\right]\\
&\leq2\gamma
E\left[I_{\{\tau_h\wedge\sigma_{2k-1}<\infty\}}\int_{\tau_h\wedge\sigma_{2k-1}}^{\tau_h\wedge(\sigma_{2k-1}+T)}
\|x(s)\|^{p-1}\|x(\tau_h\wedge\sigma_{2k-1})\|ds\right]\\
&\leq2\gamma Th^p:=C_2T^{\frac{1}{2}}.
\end{split}
\end{equation}
Next it could be deduced from \eqref{eq0} that
\begin{equation}\label{21}
\begin{split}
2&E\left[I_{\{\tau_h\wedge\sigma_{2k-1}<\infty\}}\int_{\tau_h\wedge\sigma_{2k-1}}^{\tau_h\wedge(\sigma_{2k-1}+T)}
|\langle
f(s,x(s),x(s-\tau)),x(s)-x(\tau_h\wedge\sigma_{2k-1})\rangle_H|
ds\right]\\
&\leq2E\left[I_{\{\tau_h\wedge\sigma_{2k-1}<\infty\}}\int_{\tau_h\wedge\sigma_{2k-1}}^{\tau_h\wedge(\sigma_{2k-1}+T)}
\| f(s,x(s),x(s-\tau))\|_H\|x(s)-x(\tau_h\wedge\sigma_{2k-1})\|_H
ds\right]\\
&\leq2(M+L_h)\left(1+2h/\beta\right)^2T:=C_3T^{\frac{1}{2}},
\end{split}
\end{equation}
and
\begin{equation}\label{22}
\begin{split}
2&E\left[I_{\{\tau_h\wedge\sigma_{2k-1}<\infty\}}\int_{\tau_h\wedge\sigma_{2k-1}}^{\tau_h\wedge(\sigma_{2k-1}+T)}\|g(s,x(s),x(s-\tau))\|_2^2ds\right]\\
&\leq2(M+L_h)^2\left(1+2h/\beta\right)^2T:=C_4T^{\frac{1}{2}}.
\end{split}
\end{equation}
Furthermore, taking into account Burhold-Davis-Gundy inequality,
together with  \eqref{eq1},
\begin{equation}\label{23}
\begin{split}
2E&\left[I_{\{\tau_h\wedge\sigma_{2k-1}<\infty\}}\sup\limits_{0\leq
t\leq
T}\left|\int_{\tau_h\wedge\sigma_{2k-1}}^{\tau_h\wedge(\sigma_{2k-1}+t)}\langle
x(s)-x(\tau_h\wedge\sigma_{2k-1}),
g(s,x(s),x(s-\tau))dB(s)\rangle_H\right|\right]\\
&\leq
6E\left[I_{\{\tau_h\wedge\sigma_{2k-1}<\infty\}}\int_{\tau_h\wedge\sigma_{2k-1}}^{\tau_h\wedge(\sigma_{2k-1}+T)}\|
x(s)-x(\tau_h\wedge\sigma_{2k-1})\|_H^2
\|g(s,x(s),x(s-\tau))\|_2^2ds\right]^{\frac{1}{2}}\\
&\leq6(M+L_h)(1+2h/\beta)^2T^{\frac{1}{2}}:=C_5T^{\frac{1}{2}}.
\end{split}
\end{equation}
Hence, putting \eqref{19}-\eqref{23} into \eqref{18}, it holds
\begin{equation}\label{24}
\begin{split}
E&\left[I_{\{\tau_h\wedge\sigma_{2k-1}<\infty\}}\sup\limits_{0\leq
t\leq T}
\|x(\tau_h\wedge(\sigma_{2k-1}+t))-x(\tau_h\wedge\sigma_{2k-1})\|_H^2\right]\\
& \leq (C_1+C_2+C_3+C_4+C_5)T^{\frac{1}{2}}.
\end{split}
\end{equation}
Since $w(x)$ is continuous in $V$, it must be uniformly continuous
in the closed ball $\bar{S}_h=\{x\in V: \|x\|\leq h\}$. We can
therefore choose $\delta=\delta(\epsilon)>0$ so small that
\begin{equation}\label{25}
|w(x)-w(y)|\leq\epsilon
\end{equation}
whenever $\|x-y\|\leq\delta$ with $x,y\in\bar{S}_h$.
We further choose $T =T(\epsilon,\delta,h)>0$ sufficiently small for
\begin{equation*}
(C_1+C_2+C_3+C_4+C_5)T^{\frac{1}{2}}/\delta<\epsilon.
\end{equation*}
By the Chebyshev inequality, \eqref{24} gives
\begin{equation*}
P\left(\{\tau_h\wedge\sigma_{2k-1}<\infty\}\bigcap\left\{\sup\limits_{0\leq
t\leq T}
\|x(\tau_h\wedge(\sigma_{2k-1}+t))-x(\tau_h\wedge\sigma_{2k-1})\|_H^2\geq\delta\right\}\right)<\epsilon.
\end{equation*}
Consequently,
\begin{equation*}
\begin{split}
P&\left(\{\sigma_{2k-1}<\infty,\tau_h=\infty\}\bigcap\left\{\sup\limits_{0\leq
t\leq T}
\|x(\sigma_{2k-1}+t)-x(\sigma_{2k-1})\|_H^2\geq\delta\right\}\right)\\
&=P\left(\{\tau_h\wedge\sigma_{2k-1}<\infty,\tau_h=\infty\}\bigcap\left\{\sup\limits_{0\leq
t\leq T}
\|x(\tau_h\wedge(\sigma_{2k-1}+t))-x(\sigma_{2k-1})\|_H^2\geq\delta\right\}\right)\\
&\leq
P\left(\{\tau_h\wedge\sigma_{2k-1}<\infty\}\bigcap\left\{\sup\limits_{0\leq
t\leq T}
\|x(\tau_h\wedge(\sigma_{2k-1}+t))-x(\sigma_{2k-1})\|_H^2\geq\delta\right\}\right)\\
&\leq\epsilon.
\end{split}
\end{equation*}
Recalling \eqref{17} and \eqref{27}, we further compute
\begin{equation*}
\begin{split}
P&\left(\{\sigma_{2k-1}<\infty,\tau_h=\infty\}\bigcap\left\{\sup\limits_{0\leq
t\leq T}
\|x(\sigma_{2k-1}+t)-x(\sigma_{2k-1})\|_H^2<\delta\right\}\right)\\
&=P\left(\{\sigma_{2k-1}<\infty,\tau_h=\infty\}\right)\\
&-P\left(\{\sigma_{2k-1}<\infty,\tau_h=\infty\}\bigcap\left\{\sup\limits_{0\leq
t\leq T}
\|x(\sigma_{2k-1}+t)-x(\sigma_{2k-1})\|_H^2\geq\delta\right\}\right)\\
&\geq2\epsilon-\epsilon=\epsilon.
\end{split}
\end{equation*}
Using \eqref{25}, we derive that
\begin{equation}\label{30}
\begin{split}
P&\left(\{\sigma_{2k-1}<\infty,\tau_h=\infty\}\bigcap\left\{\sup\limits_{0\leq
t\leq T}
|w(x(\sigma_{2k-1}+t))-w(x(\sigma_{2k-1}))|<\epsilon\right\}\right)\\
&\geq
P\left(\{\sigma_{2k-1}<\infty,\tau_h=\infty\}\bigcap\left\{\sup\limits_{0\leq
t\leq T}
\|x(\sigma_{2k-1}+t)-x(\sigma_{2k-1})\|<\delta\right\}\right)\\
&\geq\epsilon.
\end{split}
\end{equation}
Set
\begin{equation*}
\bar{\Omega}_k=\left\{\sup\limits_{0\leq t\leq T}
|w(x(\sigma_{2k-1}+t))-w(x(\sigma_{2k-1}))|<\epsilon\right\}.
\end{equation*}
Noting
\begin{equation*}
\sigma_{2k}(\omega)-\sigma_{2k-1}(\omega)\geq T
\end{equation*}
whenever
$\omega\in\{\sigma_{2k-1}<\infty,\tau_h=\infty\}\bigcap\bar{\Omega}_k$,
we derive from \eqref{31} and \eqref{30} that
\begin{equation*}
\begin{split}
\infty&>\epsilon\sum\limits_{k=1}^{\infty}E[I_{\{\sigma_{2k-1}<\infty,\tau_h=\infty\}}(\sigma_{2k}-\sigma_{2k-1})]\\
&\geq\epsilon\sum\limits_{k=1}^{\infty}E[I_{\{\sigma_{2k-1}<\infty,\tau_h=\infty\}\bigcap\bar{\Omega}_k}(\sigma_{2k}-\sigma_{2k-1})]\\
&\geq\epsilon T\sum\limits_{k=1}^{\infty}\epsilon=\infty
\end{split}
\end{equation*}
which is a contradiction. So \eqref{12} must hold. We observe from
\eqref{12} and (\eqref{15} that there is an $\Omega_0\subset\Omega$
with $P(\Omega_0)=1$ such that for all $\omega\in\Omega_0$
\begin{equation}\label{33}
\lim\limits_{t\rightarrow\infty}w(x(t),\omega)=0\ \ \mbox{ and} \ \
\sup\limits_{0\leq\leq t<\infty}\|x(t,\omega)\|<\infty.
\end{equation}
We shall now show that for any $\omega\in\Omega_0$
\begin{equation}\label{34}
\lim\limits_{t\rightarrow\infty}\|x(t,\omega)\|=0.
\end{equation}
If this is false, then there is some $\bar{\omega}\in\Omega_0$ such
that
\begin{equation*}
\lim\limits_{t\rightarrow\infty}\sup\|x(t,\bar{\omega})\|>0.
\end{equation*}
Whence there is a subsequence $\{x(t_k, \bar{\omega})\}_{k\geq1}$ of
$\{x(t,\bar{\omega})\}_{t\geq0}$ such that
\begin{equation*}
\|x(t_k ,\bar{\omega})\|\geq\rho,\ \ \ \ k\geq1
\end{equation*}
for some $\rho>0$. Since $\{x(t_k ,\bar{\omega})\}_{k\geq1}$ is
bounded so there must be an increasing subsequence
$\{\bar{t}_k\}_{k\geq1}$ such that $\{x(\bar{t}_k,
\omega)\}_{k\geq1}$ converges to some $z \in V$ with
$\|z\|\geq\rho$. Hence
\begin{equation*}
w(z)=\lim\limits_{k\rightarrow\infty}w(x(\bar{t}_k,\omega)).
\end{equation*}
However, by \eqref{33}, $w(z)=0$. This is a contradiction and hence
\eqref{34} must hold. Therefore,
\begin{equation*}
P\left(\lim\limits_{t\rightarrow\infty}\|x(t)\|_H=0\right)=1.
\end{equation*}
by using \eqref{eq1}. That is, the solution $x(t)$ of \eqref{eq2} is
almost surely asymptotically stable, and the proof is therefore
complete.
\end{proof}

\begin{remark}

\end{remark}

Now one example is constructed to illustrate our theory.
\begin{example}
Consider the following semilinear stochastic partial differential
equation with delay:
\begin{equation}\label{eq6}
\begin{cases}
dy(t,x)&=\frac{\partial^2}{\partial
x^2}y(t,x)dt-(y^3(t,x)+y(t,x))dt+y(t-\tau,x)\sin t dB(t),\ \
t\geq0,\ \
x\in(0,\pi),\\
y(t,x)&=\phi(t,x),\ \ 0\leq x\leq\pi,\ \ t\in[-\tau,0];\ \
y(t,0)=y(t,\pi)=0, \ \ t\geq0,
\end{cases}
\end{equation}
where $\phi\in C^2([0,\pi]\times[-\tau,0];R)$, $\tau>0$ is a
positive constant, and $B(t),t\geq0$, is a real standard Borwnian
motion.
\end{example}
We can set this problem in our formulation by taking
$H=L^2([0,\pi]), V=H_0^1([0,\pi]), V^*=H^{-1}([0,\pi]), K=R,
A(t,u)=\frac{\partial^2}{\partial x^2}u(x),f(t,u,v)=-[u^3(x)+u(x)]$
and $g(t,u,v)=v(x)\sin t$. Furthermore, the norms in $H$ and $V$ are
defined as
$\|\xi\|_H=\left(\int_0^{\pi}\xi^2(s)ds\right)^{\frac{1}{2}}$ for
$\xi\in H$ and $\|\xi\|=\left(\int_0^{\pi}\left(\frac{\partial
\xi}{\partial s}\right)^2ds\right)^{\frac{1}{2}}$ for $\xi\in V$,
respectively.
 Setting
$U(t,x)=\|x\|_H^2$ and recalling the definition of diffusion
operator $\mathcal {L}V$, it follows easily that
\begin{equation*}
\begin{split}
\mathcal {L}U(t,x,y)&=2\langle A(t,x),x\rangle+2\langle
f(t,x,y),y\rangle_H+\|g(t,x,y)\|_H^2\\
&=2\langle A(t,x),x\rangle+2\langle -x^3-x,x\rangle_H+\|y\sin t
\|_H^2\\
& \leq-2\|x\|^2-2(\|x\|^4_H+\|x\|^2_H)+\|y\|^2_H\\
&\leq-2(\|x\|^4_H+2\|x\|^2_H)+\|y\|^2_H.
\end{split}
\end{equation*}
Taking $\gamma(t)=0, w_1(x)=2(\|x\|_H^4+2\|x\|_H^2)$,
$w_2(x)=\|x\|_H^2$ and then applying Theorem \ref{Lasalle}, the
solution of \eqref{eq6} is almost surely asymptotically stable.

Now, we further take into account another Khasminskii-type condition
to give a powerful criterion  for exponential stability of
stochastic evolution delay equation, especially for highly nonlinear
cases.
\begin{theorem}\label{exponential}
Let Assumption \ref{as1} and Assumption \ref{as2} hold. Assume that
there are functions $U\in C^{1, 2}(R_+\times H;R_+)$, $W_1\in
C(R_+\times H;R_+)$, $\gamma(t), t\in R_+$, nonnegative continuous
function, and constants
$\beta_1>0,\beta_2>0,\alpha_1>\alpha_2\geq0,\alpha_3>\alpha_4>0,\mu>0$
such that

\begin{equation}\label{eq17}
\beta\|x\|_H^2\leq U(t,x)\leq\beta_2\|x\|_H^2,\ \ \ \ \ \ \ \forall
(t,x)\in R_+\times V
\end{equation}
and
\begin{equation}\label{eq18}
\mathcal {L}U(t,x,y)\leq
\gamma(t)-\alpha_1U(t,x)+\alpha_2U(t-\tau,y)-\alpha_3W_1(t,x)+\alpha_4W_1(t-\tau,y),
\ \ \ (t,x,y)\in R_+\times V\times V,
\end{equation}
where $\gamma(t)$ satisfies $\int_0^\infty\gamma(t)e^{\mu t}dt<\infty.$ Then
\begin{equation}\label{eq19}
\lim\sup\limits_{t\rightarrow\infty}\frac{1}{t}\log(E\|x(t)\|_H^2)\leq-(\mu\wedge\epsilon),
\end{equation}
where $\epsilon=\epsilon_1\wedge\epsilon_2$ while $\epsilon_1>0$ and
$\epsilon_1>0$ are the unique roots to the following equations
\begin{equation}\label{eq22}
\alpha_1=\epsilon_1+\alpha_2e^{\epsilon_1\tau}\ \ \ \ and\ \ \
\alpha_3=\alpha_4e^{\epsilon_2\tau}.
\end{equation}
In other words, the global strong solution of \eqref{eq2} is mean
square exponential stability and the Lyapunov exponent should not be
greater than $-(\mu\wedge\epsilon)$.
\end{theorem}

\begin{proof}
Noting that \eqref{eq17} and \eqref{eq18} imply \eqref{eq7} and
\eqref{eq8}, respectively, therefore, \eqref{eq2} has a unique
global strong solution on $t\geq0$. To show the desired assertion
\eqref{eq19}, compute by It\^o's  formula and \eqref{eq18} that for
$t\geq0$
\begin{equation}\label{eq20}
\begin{split}
E&(e^{\epsilon t}U(t,x(t)))\\
&=EU(0,x(0))+\epsilon E\int_0^te^{\epsilon
s}U(s,x(s))ds+E\int_0^te^{\epsilon s}\mathcal
{L}U(s,x(s),x(s-\tau))ds\\
&\leq EU(0,x(0))+\int_0^t\gamma(s)e^{\epsilon s}ds+\epsilon
E\int_0^te^{\epsilon s}U(s,x(s))ds-\alpha_1E\int_0^te^{\epsilon
s}U(s,x(s))ds\\
&+\alpha_2E\int_0^te^{\epsilon
s}U(s-\tau,x(s-\tau))ds-\alpha_3E\int_0^te^{\epsilon
s}W_1(s,x(s))ds+\alpha_4E\int_0^te^{\epsilon
s}W_1(s-\tau,x(s-\tau))ds.
\end{split}
\end{equation}
Observing that
\begin{equation*}
\begin{split}
\alpha_2\int_0^te^{\epsilon
s}U(s-\tau,x(s-\tau))ds&=\alpha_2\int_{-\tau}^{t-\tau}e^{\epsilon
s}U(s,x(s))ds\\
&\leq\alpha_2\int_{-\tau}^0e^{\epsilon
s}U(s,x(s))ds+\alpha_2e^{\epsilon\tau}\int_0^te^{\epsilon
s}U(s,x(s))ds.
\end{split}
\end{equation*}
Similarly,
\begin{equation*}
\begin{split}
\alpha_4\int_0^te^{\epsilon
s}W_1(s-\tau,x(s-\tau))ds&=\alpha_4\int_{-\tau}^{t-\tau}e^{\epsilon
s}W_1(s,x(s))ds\\
&\leq\alpha_4\int_{-\tau}^0e^{\epsilon
s}W_1(s,x(s))ds+\alpha_4e^{\epsilon\tau}\int_0^te^{\epsilon
s}W_1(s,x(s))ds.
\end{split}
\end{equation*}
Hence, in \eqref{eq20}
\begin{equation}\label{eq21}
\begin{split}
E(e^{\epsilon t}U(t,x(t)))&\leq
C_3+\int_0^t\gamma(s)e^{[\epsilon+\mu-(\mu\wedge\epsilon)s]}ds-(\alpha_1-\alpha_2e^{\epsilon\tau}-\epsilon)
E\int_0^te^{\epsilon s}U(s,x(s))ds\\
 &-(\alpha_3-\alpha_4e^{\epsilon\tau})E\int_0^te^{\epsilon
s}W_1(s,x(s))ds,
\end{split}
\end{equation}
where
\begin{equation*}
C_3=EU(0,x(0))+\alpha_2e^{\epsilon\tau}E\int_{-\tau}^0e^{\epsilon
s}U(s,x(s))ds+\alpha_4e^{\epsilon\tau}E\int_{-\tau}^0e^{\epsilon
s}W_1(s,x(s))ds.
\end{equation*}
Furthermore, by \eqref{eq22}
\begin{equation*}
E(e^{\epsilon t}U(t,x(t)))\leq
C_3+C_4e^{[\epsilon-(\mu\wedge\epsilon)]t},
\end{equation*}
where $C_4=\int_0^t\gamma(s)e^{\mu s}ds<\infty$. This certainly
implies that
\begin{equation*}
E(U(t,x(t)))\leq e^{-\epsilon t}C_3+C_4e^{-(\mu\wedge\epsilon)t}\leq
(C_3\vee C_3)e^{-(\mu\wedge\epsilon)t},
\end{equation*}
and then the desired assertion follows from \eqref{eq17}.
\end{proof}

\begin{example}
 Consider
the following stochastic delay differential equation:
\begin{equation}\label{eq24}
\begin{cases}
dy(t,x)&=\frac{\partial}{\partial x}\left(a(t,x)\frac{\partial
y(t,x)}{\partial
x}\right)dt+y(t,x)(a+by(t-\tau,x)-y^2(t,x))dt\\
&+cy(t,x)y(t-\tau,x) dB(t),\ \ t\geq0,\ \
x\in(0,\pi),\\
y(t,x)&=\phi(t,x),\ \ 0\leq x\leq\pi,\ \ t\in[-\tau,0];\ \
y(t,0)=y(t,\pi)=0, \ \ t\geq0,
\end{cases}
\end{equation}
where $\phi\in C^2([0,\pi]\times[-\tau,0];R)$, $\tau>0$ is positive
constant, $B(t),t\geq0$, is a real standard Borwnian motion.
\end{example}

Indeed, define $A(t,u)=\frac{\partial}{\partial
x}\left(a(t,x)\frac{\partial u(x)}{\partial x}\right)$, where
$a(t,x)$ is measurable in $[0,\infty)\times[0,\pi]$ and satisfy
$0<\nu\leq a(t,x)\leq\alpha$ on $[0,\infty)\times[0,\pi]$ and let
$H=L^2([0,\pi]), V=H_0^1([0,\pi]), V^*=H^{-1}([0,\pi])$, with the
usual norms in the spaces $H$ and $V$ defined as
$\|\xi\|_H=\left(\int_0^{\pi}\xi^2(s)ds\right)^{\frac{1}{2}}$ for
$\xi\in H$ and $\|\xi\|=\left(\int_0^{\pi}\left(\frac{\partial
\xi}{\partial s}\right)^2ds\right)^{\frac{1}{2}}$ for $\xi\in V$,
respectively. Setting $U(t,x)=\|x\|_H^2$, then it could be derived
that
\begin{equation*}
\begin{split}
\mathcal {L}U(t,x,y)&=2\langle A(t,x),x\rangle+2\langle
f(t,x,y),x\rangle_H+\|g(t,x,y)\|_H^2\\
&=2\langle A(t,x),x\rangle+2\langle x(a+by-x^2),x\rangle_H+\|cxy
\|_H^2\\
& \leq-2\nu\|x\|^2+2a\|x\|^2_H+\frac{1}{2}\|x\|^4_H+2b^2\|y\|^2_H-2\|x\|^4_H+\frac{1}{2}\|x\|^4_H+\frac{1}{2}c^4\|y\|^4_H\\
&\leq-2(\nu-a)\|x\|^2_H+2b^2\|y\|^2_H-\|x\|^4_H+\frac{1}{2}c^4\|y\|^4_H.
\end{split}
\end{equation*}
Hence, if $\nu-a>b^2>0$ and $c^4<2$, by Theorem \ref{exponential},
the global strong solution to \eqref{24} is exponential stability.

\end{document}